\documentclass[10pt]{article}
\usepackage{graphicx}
\usepackage{amssymb}
\usepackage{epstopdf}
\usepackage{amsmath}
\usepackage{amsfonts}
\newcommand{\braket}[2]{\langle #1,#2 \rangle}

\def\phi{{\varphi}}
\def\H{{\mathcal H}}
\def\P{{\mathcal P}}

\def\phi{{\varphi}}
\def\H{{\mathcal H}}
\def\P{{\mathcal P}}

\DeclareSymbolFont{AMSb}{U}{msb}{m}{n}
\DeclareMathSymbol{\N}{\mathbin}{AMSb}{"4E}
\DeclareMathSymbol{\Z}{\mathbin}{AMSb}{"5A}
\DeclareMathSymbol{\R}{\mathbin}{AMSb}{"52}
\DeclareMathSymbol{\Q}{\mathbin}{AMSb}{"51}
\DeclareMathSymbol{\I}{\mathbin}{AMSb}{"49}
\DeclareMathSymbol{\C}{\mathbin}{AMSb}{"43}

\def\be{\begin{equation}}
\def\ee{\end{equation}}
\def\ber{\begin{eqnarray}}
\def\eer{\end{eqnarray}}

\def\beq{\begin{equation}}
\def\eeq{\end{equation}}

\font\tenBbb=msbm10   
 \font\sevenBbb=msbm7
   \font\fiveBbb=msbm5
\newfam\Bbbfam
\textfont\Bbbfam=\tenBbb \scriptfont\Bbbfam=\sevenBbb
\scriptscriptfont\Bbbfam=\fiveBbb
\def\Bbb{\fam\Bbbfam\tenBbb}
\def\R{{\Bbb R}}

\def\d#1{{\rm d}\hbox{$\mskip 0.5mu$}#1} 
\def\ssp{\mskip 1.5mu} 
\def\Id{\ssp{\rm Id}\ssp} 
\def\fg{\leavevmode
 \raise.3ex\hbox{$\mkern4.5mu\scriptscriptstyle\rangle\!\rangle$}}

\def\og{\leavevmode
 \raise.3ex\hbox{${\scriptscriptstyle\langle\!\langle}\mkern4.5mu$}}
\def\up#1{\raise 1.5pt\hbox{#1}}  

\begin{document}

\addtolength{\textheight}{0 cm} \addtolength{\hoffset}{0 cm}
\addtolength{\textwidth}{0 cm} \addtolength{\voffset}{0 cm}

\newenvironment{acknowledgement}{\noindent\textbf{Acknowledgement.}\em}{}

\setcounter{secnumdepth}{5}
 \newtheorem{proposition}{Proposition}[section]
\newtheorem{theorem}{Theorem}[section]
\newtheorem{lemma}[theorem]{Lemma}
\newtheorem{coro}[theorem]{Corollary}
\newtheorem{remark}[theorem]{Remark}
\newtheorem{claim}[theorem]{Claim}
\newtheorem{conj}[theorem]{Conjecture}
\newtheorem{definition}[theorem]{Definition}
\newtheorem{application}{Application}

\newtheorem{corollary}[theorem]{Corollary}

\title{Remarks on multi-marginal symmetric Monge-Kantorovich problems}
\author{Nassif  Ghoussoub\thanks{Partially supported by a grant
from the Natural Sciences and Engineering Research Council of Canada.}\\
{\it\small Department of Mathematics}\\
{\it\small  University of British Columbia}\\
{\it\small Vancouver BC Canada V6T 1Z2}\\
{\it\small nassif@math.ubc.ca}\vspace{1mm}\and
Bernard Maurey\thanks{Research supported by CNAV93.}\hspace{2mm}\\
{\it\small Institut de Math\'ematiques, UMR 7586 -
CNRS}\\
{\it\small 
Universit\'e Paris Diderot - Paris 7}\\
{\it\small Paris, France}\\
{\it\small maurey@math.jussieu.fr}\\
%
}
\maketitle

\begin{abstract} Symmetric Monge-Kantorovich transport problems involving a cost function given by a family of vector fields were used by Ghoussoub-Moameni to establish polar decompositions of such vector fields into $m$-cyclically monotone maps composed with measure preserving $m$-involutions ($m\geq 2$). In this note, we  relate these symmetric transport problems to the Brenier solutions of the Monge and Monge-Kantorovich problem, as well as to the Gangbo-\'Swi\c{e}ch solutions of their multi-marginal counterparts, both of which involving quadratic cost functions.

\end{abstract}

\section{Introduction} Given Borel probability measures $\mu_i$,  $i=0,1,,...,m-1$ on domains $\Omega_i \subset \mathbb{R}^d$, and a cost function $c:\Omega_0 \times \Omega_1\times.... \times \Omega_{m-1}\rightarrow \mathbb{R}$, the multi-marginal version of Monge's optimal transportation problem is to minimize: 
\begin{equation*}
C(T_1,...,T_{m-1}):=\int_{\Omega_0}c(x_0,T_1(x_0),T_2(x_0),...,T_{m-1}(x_0))d\mu_0 \tag{\textbf{M}}
\end{equation*}
among all $(m-1)$-tuples of measurable maps $(T_1,T_2,...,T_{m-1})$, where $T_i: \Omega_0 \rightarrow \Omega_i$ pushes $\mu_0$ forward to $\mu_i$ for all $i=1,...,m-1$.
The Kantorovich formulation of the problem is to minimize:
\begin{equation*}
C(\theta):=\int_{\Omega_0 \times \Omega_1\times ... \times \Omega_{m-1}}c(x_0,x_1,x_2,...,x_{m-1})d\theta \tag{\textbf{K}}
\end{equation*}
among all probability measures $\theta$ on $\Omega_0 \times \Omega_1 \times ...\times \Omega_{m-1}$ such that the canonical projection
\begin{equation*}
\pi_i:\Omega_0 \times \Omega_1 \times...\times \Omega_{m-1} \rightarrow \Omega_i
\end{equation*}
pushes $\theta$ forward to $\mu_i$ for all $i$.  

Note that for any $(m-1)$-tuple $(T_1,T_2,...,T_{m-1})$ such that $T_{i\#}\mu_0=\mu_i$ for all $i=1, 2,...,m-1$, we can define the measure $\theta=(I,T_1,T_2,...,T_{m-1})_{\#}\mu_0$ on $\Omega_0 \times \Omega_1 \times...\times \Omega_{m-1}$, where $I:\Omega_0 \rightarrow \Omega_0$ is the identity map.  Then $\theta$ projects to $\mu_i$ for all $i$ and $C(T_1,T_2,...,T_{m-1})=C(\theta)$. In other words, (\textbf{K}) can be interpreted as a relaxed version of (\textbf{M}).

 Standard results for fairly general cost functions $c$ show that there exists a probability measure $\bar \theta$ on $\Omega_0 \times \Omega_1 \times ...\times \Omega_{m-1}$ with marginals $\mu_i$,  $i=0,1,,...,m-1$, where the supremum in (\textbf{K}) is attained. The natural question here is the following:\\

{\bf Problem (1):} For which cost functions $c$, problem $(\textbf{K})$ admits a solution $\bar \theta$ 
that is  supported 
 on a ``graph", that is a measure of the form $\bar \theta=(I,T_1,T_2,...,T_{m-1})_{\#}\mu_0$ for a suitable family of point transformations $(T_1,T_2,...,T_{m-1})$.\\

Whereas the case when $m=2$ is already well understood, the problem when $m \geq 3$ remains elusive since  existence and uniqueness in (\textbf{M}) as well as uniqueness in (\textbf{K}) are still largely open for general cost functions. There is however one important case where this problem has been resolved by Gangbo and \'Swi\c{e}ch \cite{GS}. This is when the cost function is given by 
\begin{equation}\label{quadratic}
c(x_0, x_1,..., x_{m-1})=\sum_{i=0}^{m-1}\sum_{j=i+1}^{m-1}|x_i-x_j|^2.
\end{equation}
Other extensions were also given by Pass  \cite{P5, P4, P1, P} and by Carlier-Nazaret \cite{CN}. 

In this note, we are interested in the symmetric versions of the Monge-Kantorovich problem. They are of the following type.  
 \begin{equation}\label{gm2}
{\bf K}_{\rm sym}=\sup\left\{\int_{\Omega^m}c(x_0,x_1,..., x_{m-1}) d\theta;\, \theta\in \P_{sym}(\Omega^m, \mu)\right\}
\end{equation}
 where $c$ is an appropriate cost function and $\P_{sym}(\Omega^m, \mu)$ denotes the set of all probability measures on $\Omega^m$, which are invariant under the cyclical permutation 
 \[
 \sigma (x_0, x_1,..., x_{m-1})=(x_1, x_2,..., x_{m-1}, x_0)
 \] 
 and whose marginals are all equal to a given measure $\mu$. Note that one can then assume that the cost function $c$ is cyclically symmetric, since one can replace it by its symmetrization 
 \begin{equation}
 \tilde c({\bf x})=\frac{1}{m}\sum_{i=0}^{m-1} c(\sigma^i({\bf x})), 
 \end{equation}
 and in this case, one can minimize over the set $\P(\Omega^m, \mu)$ of all probability measures on $\Omega^m$ whose all marginals are equal to $\mu$.
 
 Standard results for fairly general cost functions $c$ show that there exists $ \bar \theta \in \P_{sym}(\Omega^m,\, \mu)$, where the supremum above is attained. The natural question here is the following:\\
 
 {\bf Problem (2):} For which cost functions $c$, problem $({\bf K}_{\rm sym})$ admits as a solution a probability measure $\bar \theta$  of the form $\bar\theta=(I, S, S^2,..., S^{m-1})_\#\mu$, where $S$ is a $\mu$-measure preserving transformation on $\Omega$ such that $S^m=I$ a.e.  \\
 
 Problem (2) was resolved by Ghoussoub and Moameni for $m=2$ in \cite{GM} and for $m\geq 3$ in \cite{GM2} in the case where the cost function is of the form 
 \begin{equation}\label{mycost}
 c(x_0, x_1,..., x_{m-1})=\langle u_1(x_0), x_1\rangle +....+ \langle u_{m-1}(x_0), x_{m-1}\rangle, 
 \end{equation}
 where $u_1,..., u_{m-1}$ are bounded vector fields from $\Omega  \to \mathbb{R}^d$. Their work was in the context of establishing polar decompositions of vector fields in terms of monotone operators, which we will briefly describe below. The raison-d'\^etre of this paper is however to make a link between the results of Gangbo and \'Swi\c{e}ch dealing with the quadratic cost 
 (\ref{quadratic}) but for marginals of the form $\mu_i=\sigma^{i}_\# \mu$ for $i=0, ...., {m-1}$, and the symmetric Monge-Kantorovich problems considered by Ghoussoub and Moameni for the cost (\ref{mycost}).

\subsection*{Polar decompositions}

Recall that a vector field $u:\Omega \to \mathbb{R}^d$ on a domain $\Omega$ in $\mathbb{R}^d$ is said to be monotone on $\Omega$  if for all $(x,y)$ in $\Omega$, 
\begin{equation}\label{mon}
\hbox{$\langle x-y, u(x) - u(y)\rangle \geq 0$.}
\end{equation}
A result of E. Krauss  \cite{Kra} states that  a map $u:\Omega \to \mathbb{R}^d$ is  {\it   monotone}  if and only if
  \begin{equation}
  \hbox{$u(x)=\nabla_2H(x, x)$ for all $x\in \Omega$,}
  \end{equation}
   where $H$ is a  concave-convex anti-symmetric Hamiltonian on $\mathbb{R}^d\times \mathbb{R}^d$. More recently, Galichon-Ghoussoub \cite{GG} extended Krauss' result to the case of {\it $m$-cyclically monotone vector fields}, where $m$ is a fixed integer larger than $2$. Recall that these are the maps $u$ from $\Omega$ to $\mathbb{R}^d$ that satisfy for any $m+1$ points $(x_i)_{i=0}^m$  in $\Omega$ with $x_0=x_m$, the inequality
\begin{equation}\label{mcm}
\hbox{$\sum\limits_{k=0}^{m-1}\langle u(x_{k+1}), x_{k+1}-x_{k} \rangle \geq 0$.}
\end{equation}
For that, Galichon and Ghoussoub  consider the class ${\mathcal H}_m(\Omega)$ of all {\em $m$-cyclically antisymmetric Hamiltonians} on $\Omega^m$, that is the set
\begin{equation*}
\hbox{${\mathcal H}_m(\Omega)=\{H\in C(\Omega^m; \mathbb{R});\, \sum_{i=1}^{m}H\left( \sigma^{i-1}\left( {\bf x}\right) \right) =0\, $ for all ${\bf x}\in \Omega^m\}$,}
\end{equation*}
where $\sigma$ is the cyclical permutation $\sigma (x_0, x_1,..., x_{m-1})=(x_1, x_2,..., x_{m-1}, x_0)$. 
They then show that if a vector field $u$ is {\it $m$-cyclically monotone}, then there exists a Hamiltonian $H\in {\mathcal H}_m(\Omega)$ such that 
\begin{equation}
\hbox{$u(x)=\nabla_m H(x,x,...,x)$ for all $x\in \Omega$.}
\end{equation}
Moreover, $H$ can be assumed to be concave in the first variable, convex in the last $(m-1)$ variables,  though only {\it $m$-cyclically sub-antisymmetric} on $\Omega^m$, that is $\sum_{i=1}^{m}H (\sigma^{i-1}( {\bf x})) \leq 0$ for all ${\bf x}\in \Omega^m.$

It is worth comparing the above to a classical theorem of Rockafellar \cite{Ph}, which yields that
 a single-valued map $u$ from $\Omega$ to $\mathbb{R}^d$  is  a {\it maximal cyclically monotone operator} (i.e., satisfies (\ref{mcm}) for every $m\geq 2$), if and only if
\begin{equation}
\hbox{$u (x)=\nabla \phi (x)$ on $\Omega$, where $\phi:\mathbb{R}^d \to \mathbb{R}$ is a convex function. }
\end{equation}

More remarkable is the polar decomposition that Y. Brenier \cite{Br} establishes for a general  non-degenerate vector field, and which follows from his celebrated mass transport theorem. Recall that a mapping $u: \Omega \rightarrow \mathbb{R}^d$ is said to be {\it non-degenerate} if  the inverse image $u^{-1}(N)$ of  every zero-measure $N\subseteq \mathbb{R}^d$ has also zero measure. Brenier proved that any non-degenerate vector field $u \in L^{\infty}(\Omega, \mathbb{R}^d)$ can be decomposed  as
\begin{equation}\label{pol.0}
\hbox{$u(x)=\nabla \phi \circ S(x)$ \,\, a.e. in $\Omega$,}
\end{equation}
with
$\phi:\mathbb{R}^d\rightarrow \mathbb{R}$  being a  convex function and  $S:\bar {\Omega}\rightarrow \bar {\Omega}$ a measure preserving transformation.

Later, and in  the same spirit as  Brenier's, Ghoussoub and Moameni established in \cite{GM} another decomposition for non-degenerate vector fields,  which can be seen as the general version of Krauss' characterization of monotone operators. Assuming the boundary $\partial \Omega$ has measure zero, they show that if $u \in L^{\infty}(\Omega, \mathbb{R}^d)$ is a non-degenerate vector field,  then there exists a measure preserving transformation $S:
{\Omega}\rightarrow  {\Omega}$ such that $S^2=I$ (i.e., an involution), and a globally Lipschitz  anti-symmetric concave-convex Hamiltonian $H: \mathbb{R}^d \times \mathbb{R}^d \to \mathbb{R}$ such that
\begin{equation} \label{pol.1}
u(x)= \nabla_2 H(x, Sx) \quad \text{ a.e. } \, x \in \Omega.
\end{equation}
In other words, up to a measure preserving involution, essentially every bounded vector field is monotone, where the latter correspond to when $S$ is the identity map on $\Omega$.

More recently, Ghoussoub-Moameni \cite{GM2} extended this result by showing that any family $u_1,..., u_{m-1}$ of non-degenerate bounded vector field can be represented as 
 \begin{equation}\label{pol.2}
\hbox{$u_i(x)=\nabla_{i+1}H(x, Sx, S^2x,...,S^{m-1}x) \,\, \text{ a.e. } \, x \in \Omega$\quad for $i=1,..., m-1$,}
\end{equation}
where $H\in \H_m(\Omega)$ and $S$ is a measure preserving $m$-involution (i.e., $S^m=I$). Moreover, $H$ could be replaced by a Hamiltonian that is  concave in the first variable and convex in the other $(m-1)$-variables, though only $m$-cyclically sub-antisymmetric. 

The proofs of the representations (\ref{pol.1}) (when $m=2$) and of (\ref{pol.2}) (when $m\geq 3$) rely on symmetric versions of the Monge problem and of its multi-marginal Monge-Kantorovich version as mentioned above. We shall give here another formulation that is closer to the original Monge and multi-marginal Monge-Kantorovich
 problems corresponding to a quadratic cost.

\section{Mass transport with quadratic cost in the presence of symmetry} 

Given two probability measures with finite second moment $\mu_0,\mu_1$ on $\mathbb{R}^d$, with  $X:={\rm support} (\mu_0)$ and $Y:={\rm support} (\mu_1)$, the Wasserstein distance $W_2(\mu_0, \mu_1)$ between them is defined by the formula 
\begin{equation}\label{Wass}
W_2(\mu_0, \mu_1)^2=\inf \left\{\int_{X}{|x-T(x)|^2}d\mu_0(x);\,  T\in {\cal S}(\mu_0, \mu_1) \right\}
\end{equation}
where ${\cal S}(\mu_0, \mu_1)$ is the class of all Borel measurable maps $T: X\to Y$ such that  $T_\#\mu_0=\mu_1$, i.e., those which satisfy the change of variables formula,
\begin{equation}
\int_{Y}{h(y)d\mu_1(y)}=\int_{X}{h(T(x))d\mu_0(x)},\ \ \ \text{for every}\ \ h \in C(Y).
\end{equation}
Whether the infimum describing the Wasserstein distance $W_2(\mu_0, \mu_1)$ is achieved by an optimal map $\bar T$ is a variation on the original mass transport problem of G. Monge, who inquired about finding the optimal way for rearranging $\mu_0$ into $\mu_1$ against the cost function $c(x) =|x|$. 
Our cost function here  $c(x) =\frac{1}{2}|x|^2$ is quadratic, and the existence, uniqueness and characterization of an optimal map that we give below,  was established by Y. Brenier. 

\begin{theorem} {\bf (Brenier)} \label{Brenier} Assume $\mu_0$ is absolutely continuous with respect to Lebesgue measure, 
then there exists a unique optimal map  $\bar{T}$ in ${\cal S}(\mu_0, \mu_1)$, where the infimum in (\ref{Wass}) is achieved. Moreover, the map $\bar T: X\to Y$ is one-to-one and onto, $\mu_0$ a.e., and is equal to $\nabla \phi$ $\mu_0$ a.e on $X$, for some convex function $\phi:\mathbb{R}^d\to \mathbb{R}$. 

Moreover, the Brenier map $\bar{T}$ is the unique map (up to $\mu_0$ a.e. equivalence) of the form $\nabla \phi$ with $\phi$ convex such that $\nabla \phi_ \# \mu_0 =\mu_1$. 
\end{theorem}

Now we consider the above theorem in the presence of symmetry. 

\begin{corollary} Let $\mu$ be a probability measure on $\mathbb{R}^d$ that is absolutely continuous with respect to Lebesgue measure, and let $\tilde \mu$ be its image by a self-adjoint unitary transformation $\sigma$ on $\mathbb{R}^d$ (i.e.,  $\sigma= \sigma^*$ and $ \sigma^2=I$). Then, there exists a convex function $\phi:\mathbb{R}^d\to \mathbb{R}$ such that 
\begin{equation}
\hbox{$\nabla \phi_\# \mu =\tilde \mu$ \quad and \quad $\phi^*(\sigma (x))=\phi (x)$ for $x\in \mathbb{R}^d$,}
\end{equation}
where $\phi^*$ is the Legendre transform of $\phi$.
\end{corollary}
\noindent {\bf Proof:} The above theorem yields a convex function  $\phi:\mathbb{R}^d\to \mathbb{R}$ such that $\nabla \phi_\# \mu =\tilde \mu$. Recall that if  $\phi^*$ is the Legendre transform of $\phi$, then $\nabla \phi^*=(\nabla \phi)^{-1}$. Hence the function $\psi(x)= \phi^*(\sigma (x))$, which is also convex has a gradient $\nabla \psi=\sigma ^*\circ \nabla \phi^*\circ \sigma$, which  also maps $\mu$ onto $\tilde\mu$. By the uniqueness property, we have $\nabla \psi=\nabla \phi$, which means that --up to a constant-- $\phi^*(\sigma (x))=\phi (x)$ for all $x\in \mathbb{R}^d$.\\

In the rest of this section, we try to connect the above corollary to the following polar decomposition established in \cite{GM}.

 \begin{theorem}\label{main} {\bf (Ghoussoub-Moameni)} Let $\Omega$ be an open bounded set in $\mathbb{R}^d$ such that $\partial \Omega$ has zero Lebesgue measure.
\begin{enumerate}
\item  If $u \in L^{\infty}(\Omega, \mathbb{R}^d)$ is a non-degenerate vector field,  then there exists a measure preserving transformation $S:\bar
{\Omega}\rightarrow \bar {\Omega}$ such that $S^2=I$ (i.e., an involution), and a globally Lipschitz  anti-symmetric concave-convex Hamiltonian $H: \mathbb{R}^d \times \mathbb{R}^d \to \mathbb{R}$ such that
\begin{equation} \label{rep.100}
u(x)= \nabla_2 H(x, Sx) \quad \text{ a.e. } \, x \in \Omega.
\end{equation}
The involution $S$ is obtained by solving the following variational problem
\begin{equation}
\sup\left\{\int_\Omega \langle u(x), Sx\rangle \, dx;\, S\,  \hbox{\rm is a measure preserving involution on}\,  \Omega\right\}.
\end{equation}
\item Moreover, $u$ is a  monotone map if and only if there is a representation as (\ref{rep.100}) with $S$ being the identity.
\end{enumerate}
\end{theorem}

\subsection*{Lagrangians and Hamiltonians associated to monotone maps} 

An important example of an involutive transformation is the transpose $\sigma$ on $\mathbb{R}^d\times \mathbb{R}^d$ defined by $\sigma (x, y)=(y,x)$, since the convex functions $L$ on $\mathbb{R}^d\times \mathbb{R}^d$ satisfying 
$ L^*(y,x)=L(x, y)$ for all $(x,y)\in \mathbb{R}^d\times \mathbb{R}^d$ are connected to central notions in nonlinear analysis and PDEs \cite{Gh}. 

Duality theory is at the heart of this concept and it is therefore enlightening to describe it in the case where $\mathbb{R}^d$ is replaced by any reflexive Banach space $X$. Recall from  \cite{Gh}  the notion of a vector field $\bar \partial L$ that is derived from a convex lower semi-continuous Lagrangian on phase space $L:X\times X^*\to \mathbb{R}\cup\{+\infty\}$ in the following way:  for each $x\in X$, the  --possibly empty-- subset  $\bar \partial L(x)$ of $X^*$ is defined as
 \begin{eqnarray}
\bar \partial L(x): = \{ p \in X^*;  (p, x)\in \partial L(x,p) \}.
\end{eqnarray}
 Here $\partial L$ is the subdifferential of the convex function $L$ on $X\times X^*$, which should not be confused with $\bar \partial L$. Of particular interest are those vector fields derived from  {\it self-dual Lagrangians}, i.e., those convex lower semi-continuous Lagrangians $L$ on $X\times X^*$ that satisfy the following duality property:
  \begin{equation}
  L^*(p,x)=L(x, p) \quad \hbox{\rm for all $(x,p)\in X\times X^*$},
  \end{equation}
   where here $L^*$ is the Legendre transform in both variables, i.e.,
  \[
  L^*(p,x)= \sup  \{ \braket{ y}{p }+  \braket{x}{q }-L(y, q): \,  (y,q)\in X\times X^*\}.
  \]
  Such Lagrangians satisfy the following basic property:
 \begin{equation}\label{obs.1}
\hbox{$ L(x, p)-\langle x, p\rangle\geq 0$   for every $(x, p) \in X\times X^{*}$.}
 \end{equation}
  Moreover,
 \begin{equation}\label{obs.2}
\hbox{  $L(x, p)-\langle x, p\rangle =0$ if and only if $(p, x)\in \partial L(x,p)$,}
\end{equation}
which means  that the associated vector field at $x \in X $ is simply
\begin{eqnarray}
\bar \partial L(x):= \{ p \in X^*; L(x,p)- \langle x,p \rangle=0 \}.
\end{eqnarray}
These so-called  {\it selfdual vector fields} are natural but far reaching extensions of subdifferentials of convex lower semi-continuous functions. Indeed, the most
basic selfdual Lagrangians are of the form 
\[
L(x,p)= \varphi (x)+\varphi^*(p),
\]
 where  $\varphi$ is a convex and lower semi-continuous function on $X$, and $\varphi^*$
is its Legendre conjugate on $X^*$,  in which case  
\[
\bar \partial  L(x)= \partial \varphi (x).
\]
More interesting  examples of self-dual Lagrangians are of the form 
\[
L(x,p)= \varphi (x)+\varphi^*(-\Gamma x+p),
\]  where $\varphi$  is as above,
and  $\Gamma: X\rightarrow X^*$ is a skew adjoint operator. The corresponding selfdual vector field is then
\begin{equation}
\bar \partial L(x)=\Gamma x+ \partial \varphi (x).
\end{equation}

 Actually, it turned out that any {\it maximal monotone operator} $A$ 
  is a self-dual vector field and vice-versa  \cite{Gh}. That is,  there exists a selfdual Lagrangian $L$ such that $A=\bar \partial L$. This fact was proved and reproved by several authors. See for example, R. S. Burachik and B. F. Svaiter  \cite{BS}, B. F. Svaiter \cite{S}.
      
 This result means that self-dual Lagrangians can be seen as the {\it potentials} of maximal monotone operators, in the same way as the Dirichlet integral  is the potential of the Laplacian operator (and more generally as any convex lower semi-continuous energy is a potential for its own subdifferential). Check \cite{Gh} to see how this characterization leads to variational formulations and resolutions of most equations involving monotone operators.

\begin{proposition} \label{Kra} Let $u: \Omega \to {\bf R}^d$ be a possibly set-valued map. The following properties are then equivalent:
\begin{enumerate}

\item $u$ is a maximal monotone map with domain $\Omega$.
\item There exists a convex self-dual Lagrangian on $\mathbb{R}^d\times \mathbb{R}^d$ such that 
$u(x)= \bar \partial L(x)$ for all $x\in \Omega$. In other words, 
\begin{equation}
\hbox{$(p,x)\in \partial L(x,p)$ if and only if $p\in u(x)$.}
\end{equation}
\item There exists a concave-convex anti-symmetric Hamiltonian $H$ on $\mathbb{R}^d \times \mathbb{R}^d$ such that \begin{equation}
\hbox{$u(x) \subset  \partial_2 H(x,x)$  for $x\in \Omega$.} 
\end{equation}
\end{enumerate}
\end{proposition}
\noindent {\bf Sketch of proof:}  
Assuming $u$ is maximal monotone with domain $\Omega$, we consider 
its associated Fitzpatrick function, that is 
\begin{equation}
N(p,x)=\sup\{\langle p,y\rangle +\langle q, x-y\rangle; \, (y,q)\in {\rm Graph}(u)\}.
\end{equation}
 It is known and easy to see that
 \begin{equation}
\hbox{$N^*(x,p) \geq N(p,x)\geq \langle x, p\rangle$ for every $(x,p)\in {\mathbb{R}}^d\times {\mathbb{R}}^d$,}
\end{equation}
and that 
\[
\hbox{$N(p,x)=\langle x, p\rangle$ if and only if $(x,p)\in  {\rm Graph}(u)$}.
\]
Now consider the following Lagrangian on $\mathbb{R}^d\times \mathbb{R}^d$, which interpolates between $N$ and $N^*$,
 \[
 L(p,x):=\inf \left\{\frac{1}{2}N(p_1, x_1)+\frac{1}{2}N^*(x_2, p_2)+\frac{1}{8}\|x_1-x_2\|^2+\frac{1}{8}\|p_1-p_2\|^2\right\},
\]
where the infimum is taken over all couples $(x_1, p_1)$ and $(x_2, p_2)$ such that 
 \[
(x, p)=\frac{1}{2}(x_1, p_1) + \frac{1}{2}(x_2, p_2).
 \]
It can be shown  (see \cite{Gh}, p. 92) that $L$ is a self-dual Lagrangian on $\mathbb{R}^d\times \mathbb{R}^d$,  in such a way that 
\begin{equation}
\hbox{$N^*(x,p) \geq L(p, x)\geq N(p,x)\geq \langle x, p\rangle$ for every $(x,p)\in {\mathbb{R}}^d\times {\mathbb{R}}^d$,}
\end{equation}
which means that if $(x, p) \in \Omega \times \R^d$, then $L(p,x)=\langle x, p\rangle$ if and only if $(x, p)\in {\rm Graph}(u)$, that is when $p \in u(x)$. 

To show that 2) implies 3), it suffices to take the Legendre transform of $L$ with respect to the second variable, i.e., 
\[
K_L(x,y)=\sup\{\langle y, p\rangle - L(x, p); p\in \mathbb{R}^d \}.
\]
It is clearly concave-convex. The selfduality of $L$ yields that $K_L$ is (at least) sub-antisymmetric, i.e., 
\[
\hbox{$K_L(x,y)\leq -K_L(y,x)$ for  $(x,y)\in \mathbb{R}^d\times \mathbb{R}^d$.}
\]
Note that $K_L(x,y)\leq H_L(x,y)=\frac{1}{2}(K_L(x,y)-K_L(y,x))$ is anti-symmetric, and that $u(x)\in \partial_2 H_L(x,x)$  in $\Omega$ as in the theorem of Krauss mentioned above. 

Finally, assuming 3) that is $u(x)\in \partial_2 H(x,x)$, where $H$ is convex in the second variable, we have for any $x, y \in \Omega$, any $p\in   \partial_2 H_L(x,x)$ and $q\in   \partial_2 H_L(y,y)$
\[ 
\hbox{$H(x,y) \geq H(x,x) + \langle p, y-x\rangle$ and $H(y,x) \geq H(y,y) + \langle q, x-y\rangle.$}
\]
Since $H$ is anti-symmetric (hence $H(x,x)=0$), this implies --by adding the two inequalities-- that 
\[
0 \geq \langle u(x) -u(y), y-x\rangle,
\]
which means that $u$ is monotone. \hfill $\Box$\\

Now we note that Monge transport problems provide a natural way to construct selfdual Lagrangians, hence general monotone operators.

\begin{corollary} Let $\mu$ be a probability measure on the product space $\Omega:=\Omega_1\times \Omega_2\subset\mathbb{R}^d\times \mathbb{R}^d$, and let $\tilde \mu$ be the probability measure on $\tilde \Omega:=\Omega_2\times \Omega_1$ obtained as the image of $\mu$ by the transformation $\sigma (x, y)=(y,x)$. Then, there exists a self-dual Lagrangian $L:\mathbb{R}^d\times \mathbb{R}^d \to \mathbb{R}$ such that $\nabla L_\# \mu=\tilde\mu$. 

Moreover, the Monge-Kantorovich problem 
\begin{equation}
\sup\{\int_{\Omega \times \tilde \Omega}[\langle x, y'\rangle +\langle y, x'\rangle]\, d\theta (x,y, y', x'); \theta \in {\mathcal P}(\Omega \times \tilde \Omega), \, \, \,  {\rm proj}_1\theta =\mu,  {\rm proj}_2\theta =\tilde\mu\}
\end{equation}
has a solution $\theta$ that is supported on the self-dual Lagrangian manifold 
\begin{equation}
\{((x,y), (x', y'))\in \mathbb{R}^{2d} \times \mathbb{R}^{2d}; L(x, y)+L(x', y')=\langle x, y'\rangle +\langle y, x'\rangle\}.
\end{equation}
\end{corollary}

\noindent{\bf Proof:} This is a direct application of the above corollary. The map $\nabla L$ then solves the Monge problem
\begin{equation}
\int_{\Omega_1\times \Omega_2}|(x,y)-\nabla L(x,y)|^2\, d\mu (x,y)=\inf\{\int_{\Omega_1\times \Omega_2}|(x,y)-T(x,y)|^2\, d\mu (x,y); \, T_\#\mu =\tilde\mu\}.
\end{equation}
The self-duality of $L$ follows from applying the above corollary to the transformation $\sigma (x, y)= (y,x)$. \hfill $\Box$

\subsection*{Involutions as a byproduct of mass transport between graphs} 

Suppose now that  $\mu=(I \otimes u)_\#dx$, where $dx$ is normalized Lebesgue measure on a bounded domain $\Omega$ in $\mathbb{R}^d$, and  $u:\Omega \to \mathbb{R}^d$ is a vector field in such a way that $\mu$ is  supported by the graph $G:=\{(x, u(x));\, x\in \Omega\}$ and let $\tilde \mu =\sigma_\#\mu$. Now any map $T$ pushing $\mu$ onto $\tilde \mu$ can be parameterized by an application $S:\Omega \to \Omega$ via the formula:
\begin{equation}\label{rep}
T: (x, u(x)) \to (u(Sx), Sx),
\end{equation}
and the Monge problem between $\mu$ and $\tilde \mu$ can then be formulated as 
 \begin{equation}\label{gr}
\inf\left\{\frac{1}{2}\int_\Omega [|u(Sx)-x|^2 +|u(x)-S(x)|^2] \, dx;\, S\,  \hbox{\rm measure preserving transformation on}\,  \Omega\right\}.
\end{equation}
Assume now that --just like in the non-degenerate case-- there exists a self-dual Lagrangian $L$ such that $\nabla L_\#\mu=\tilde \mu$. This means that there exists $S:\Omega \to \Omega$ such that 
\[
\hbox{$\nabla L(x, u(x))=(u(Sx), Sx)$ for a.e. $x\in \Omega$.}
\] 
Since $\nabla L^*=\sigma^*\circ  \nabla L \circ \sigma$ and $\nabla L^*(u(Sx), Sx)=(x, u(x))$ a.e., we have for a.e. $x\in \Omega$,
\[
(x, u(x))=\nabla L^*(u(Sx), Sx)=\sigma^*\circ \nabla L(Sx, u(Sx))=\sigma^*(u(S^2x), S^2x)=(S^2x, u(S^2x)).
\]
It follows that $S$ is a measure preserving involution on $\Omega$. In other words, problem (\ref{gr}) is equivalent to the problem 
\begin{equation}
\sup\left\{\int_\Omega \langle u(x), Sx\rangle \, dx;\, S\,  \hbox{\rm measure preserving involution on}\,  \Omega\right\}, 
\end{equation}
which was used by Ghoussoub-Moameni to establish the polar decomposition (\ref{pol.1}).  Actually, if one now considers the Legendre transform $H$ of $L$ with respect to the second variable, then as noted above, $H$ is sub-anti-symmetric, can be assumed to be anti-symmetric, and satisfies,
\[
\hbox{$L(x, u(x))+H(x, Sx)=\langle u(x), Sx \rangle$ for all $x\in \Omega$}.
\]
This then yields the polar decomposition $u(x)=\nabla_2H(x, Sx)$ a.e., which can be then seen as a (self-dual) mass transport problem between the measure $\mu$ supported by the graph of $u$ and its transpose. Unfortunately, the measure $\mu$ is too degenerate to fall under the framework where we have uniqueness in Brenier's theorem, hence the need to give a direct proof of the result as in \cite{GM}.

In order to link the polar decomposition with the symmetric Monge-Kantorovich problem, we note that measurable functions $S:\Omega \to \Omega$, whose graphs $\{(x, Sx); x\in \Omega\}$ are the support of measures in ${\mathcal P}^{\mu}_{{\rm sym}}(\Omega \times \Omega)$ can be characterized in the following way.

\begin{lemma} Let $S:\Omega \to \Omega$ be a measurable map, then the following are equivalent:

\begin{enumerate} 
\item The image of $\mu $ by the map $x\to (x, Sx)$ belongs to ${\mathcal P}^{\mu}_{{\rm sym}}(\Omega \times \Omega)$.

\item $S$ is $\mu$-measure preserving and $S^2x=x$ $\mu$-a.e.

\item $\int_\Omega H(Sx, x)\, d\mu (x)=0$ for every Borel measurable, bounded and antisymmetric function $H$ on $\Omega\times \Omega$. 

\end{enumerate}
\end{lemma}

\textbf{Proof.} It is clear that 1) implies 3), while 2) implies 1). We now prove that 2) and 3) are equivalent. Assuming first that $S$ is measure preserving such that $S^2=I$ a.e, then  for every anti-symmetric
$H$ in $L^\infty(\Omega \times \Omega)$, we have
\begin{eqnarray*}
\int_{\Omega} H(x,S(x)) \, d\mu (x)=\int_{\Omega} H(S(x), S^2(x)) \, d\mu (x)=\int_{\Omega} H(S(x),x) \, d\mu (x)
=-\int_{\Omega} H(x,S(x))\, d\mu (x),
\end{eqnarray*}
hence $\int_{\Omega} H(x,S(x)) \, d\mu (x)=0$. 

Conversely, if $\int_{\Omega} H(x,S(x)) \, d\mu (x)=0$ for every anti-symmetric
$H$, then it suffices to take $H(x,y)=f(x)-f(y)$, where $f$ is any continuous function on $\Omega$ to conclude that $S$ is necessarily $\mu$-measure preserving. On the other hand, if one considers the anti-symmetric functional  $$
H(x, y)=|S(x)-y|-|S(y)-x|,
$$
 then $
 0=\int_{\Omega} H(x,S(x)) \, d\mu (x)=\int |S^2(x)-x| d\mu (x),$ 
which clearly yields that $S$ is an involution $\mu$-almost everywhere. $\hfill$ $\Box$\\

We now give the following variational formulation for monotone operators, which shows that they are in some way orthogonal to involutions. Denote by ${\mathcal S}_2(\Omega, \mu)$ the set of all  $\mu$-measure preserving involutions on $\Omega$. It is an easy exercise to show that ${\mathcal S}_2(\Omega, \mu)$ is a closed subset of a sphere of $L^2(\Omega, {\mathbb{R}}^d)$.  In order to simplify the exposition, we shall assume  that $d\mu$ is Lebesgue measure $dx$ normalized to be a probability on the bounded open set $\Omega$, and $\mu$ can and will then be dropped from all notation. We shall also assume that the boundary of $\Omega$ has measure zero. 

\begin{proposition} Let $u:\Omega \to {\bf R}^d$ be a vector field in $L^2(\Omega, {\bf R}^d)$. 
The following properties are then equivalent:
\begin{enumerate}

\item $u$ is monotone a.e. on $\Omega$, that is there exists a measure zero set $N$ such that the restriction of $u$ to $\Omega \setminus N$ is monotone.

\item $\sup\big\{\int_\Omega \langle u(x), Sx-x\rangle dx; \, S\in {\mathcal S}_2(\Omega, \mu)\big\} =0$.

\item The projection of $u$ on ${\mathcal S}_2(\Omega, \mu)$ is the identity map, that is
\[
\int_\Omega |u(x)-x|^2dx=\inf\{\int_\Omega |u(x)-S(x)|^2dx; S\in {\mathcal S}_2(\Omega, \mu)\}.
\]

\item $\sup\{ \int_{\Omega \times \Omega} \langle u(x), y\rangle d\pi (x,y); \, \pi \in {\mathcal P}_{\rm sym}(\Omega \times \Omega, dx)\} =\int_\Omega \langle u(x), x\rangle \, dx .$
  \end{enumerate}
\end{proposition}
\noindent {\bf Proof:}  Assume $u$ is a monotone map and use Proposition \ref{Kra} to write it --modulo an obvious abuse in notation-- as $u (x)= \nabla_2 H(x,x)$ a.e. on $\Omega$, where $H$ is anti-symmetric and convex in the second variable. We can write for any measurable map $S$, 
\[ H(x,S(x)) \geq H(x,x) + \langle \nabla_2 H(x,x), S(x)-x)\rangle.
\]
It follows that 
$$\int_\Omega \langle u(x), x-S(x)\rangle dx=\int_\Omega \langle \nabla_2 H(x,x), x-S(x)\rangle dx \geq \int_\Omega [H(x,x)-H(x, S(x))]dx. $$
Note that $\int_\Omega H(x, S(x))dx=0$ since $S$ is measure preserving and $S^2=I$. 
Similarly, $\int_\Omega H(x,x)dx=0$ and it then follows that $\int_\Omega \langle u(x), x-Sx\rangle dx \geq 0$. Finally, by taking $Sx=x$, one can then see that the supremum in 2) is equal to zero.

Suppose now that 2) holds. In order to show 1) i.e., that $u$ is monotone, consider a pair $x_1,x_2$ in $\Omega$ and $R$ small enough so that $B(x_i,R) \subset \Omega$ for $i \in \{1,2\}$. Define the measure preserving involution $S_R$ via
$$S_R(x)=\left\{
    \begin{array}{lll}
         x-x_2+x_1 & \mbox{if } x\in B(x_2,R) \\
         x-x_1+x_2 & \mbox{if } x\in B(x_1,R) \\
         x & \mbox{otherwise}.
    \end{array}
\right. $$
Since  $\int_\Omega \langle u(x), x-S_R(x)\rangle dx \geq 0$, we have
$$\left\langle x_1 - x_2,  \int_{B_1} u(x_1 +Ry)dy - \int_{B_1} u(x_2 +Ry) dy \right\rangle \geq 0. $$
Since $u$ is Lebesgue integrable, almost every point $x \in \Omega$ is a Lebesgue point, which means that $ u(x) = \lim_{R \rightarrow 0} \mid B \mid ^{-1} \int_B u(x+Ry) dy.$
This leads to 
$\langle x_1-x_2, u(x_1)-u(x_2) \rangle \geq 0\quad \mbox{for a.e. } x_1, x_2 \in \Omega.$

By developing the square, it is clear that property 2) is equivalent to 3), which says that the identity map is the projection of $u$ on the closed subset ${\mathcal S_2}(\Omega, \mu)$ of the sphere of $L^2(\Omega, {\bf R}^d)$, that is  ${\rm dist} (u, {\cal S}_2(\Omega, \mu))=\|u-I\|_2$. In other words,  
$
\int_\Omega |u(x)-x|^2dx=\inf\{\int_\Omega |u(x)-S(x)|^2dx; S\in {\mathcal S}_2(\Omega, \mu)\}.
$

For 1) implies 4) assume $u$ is monotone and observe that for any probability $\pi$ in ${\mathcal P}(\Omega \times \Omega)$ with marginals $dx$, we have
\begin{eqnarray*}
\int_{\Omega \times \Omega}\langle u(x), y-x\rangle d\pi (x,y)&=&  \int_{\Omega \times \Omega}\langle u(x)-u(y), y-x\rangle d\pi+\int_{\Omega \times \Omega}\langle u(y), y-x\rangle d\pi\\
&\leq& \int_{\Omega \times \Omega}\langle u(y), y-x\rangle d\pi (x,y).
\end{eqnarray*}
Since $\pi$ is symmetric, we have then that
$2\int_{\Omega \times \Omega}\langle u(x), y-x\rangle d\pi (x,y) \leq 0$. The fact that the supremum is zero follows from simply taking the probability measure supported on the diagonal of $\Omega \times \Omega$. \\

Finally, note that 4) implies 2) by considering for any $S\in {\mathcal S}_2(\Omega, \mu)$ the symmetric measure $d\pi$ on $\Omega\times \Omega$ that is the image of Lebesgue measure by the map $x\to (x, Sx)$. Note that 
\begin{eqnarray*}
\int_{\Omega\times \Omega} f(x,y) d\pi (x,y)=\int_\Omega f(x, Sx) dx=\int_\Omega f(Sx, x) dx
=\int_{\Omega\times \Omega} f(x,y) d\pi (y,x).  
\end{eqnarray*}
\hfill $\Box$

Back to the case of  a general vector field, the above then shows that the variational problem used in the polar decomposition 
\begin{equation}
\sup\left\{\int_\Omega \langle u(x), Sx\rangle \, dx;\, S\,  \hbox{\rm measure preserving involution on}\,  \Omega\right\}, 
\end{equation}
is nothing but a symmetric Monge-Kantorovich problem 
\begin{equation}
\sup\{\int_{\Omega\times \Omega} \langle u (x), y\rangle  d\pi; \pi\in \P_{sym}(\Omega\times \Omega, dx) \},
\end{equation}
where the cost function is given by $c(x, y)= \langle u (x), y\rangle$. In the case where $u$ is monotone then the involution where the supremum is attained is simply the identity.

\section{Multidimensional Monge-Kantorovich Theorems}

In this section, we are interested in relating the Gangbo-\'Swi\c{e}ch solution of the multidimensional Monge-Kantorovich with quadratic cost to the following recent result of Ghoussoub-Moameni \cite{GM2}.

\begin{theorem}\label{main} \label{GM}
Given a probability measure $\mu$ on $\Omega$ and bounded vector fields $u_1, u_2,...., u_{m-1}$ from $\Omega$ to $\mathbb{R}^d$ that are $\mu$-non-degenerate, then
\begin{enumerate} 
\item The symmetric Monge-Kantorovich problem 
\begin{equation}
{\bf K_{\rm sym}}=\sup\{\int_{\Omega^m}\left[\langle u_1(x_0), x_1\rangle +...+\langle u_{m-1}(x_0), x_{m-1}\rangle \right]  d\pi; \pi\in \P_{sym}(\Omega^m,\mu) \}\label{dual}
\end{equation}
attains its maximum at a measure of the form $\bar \theta=(I,S,S^2,..,S^{m-1})_{\#}\mu$,  where $S$ is a $\mu$-measure preserving transformation on $\Omega$ such that $S^m=I$ a.e.  

\item  There exists a Hamiltonian $H\in \H_m(\Omega)$ such that for $i=1,..., m-1$, 
\begin{equation}\label{rep}
u_i(x)=\nabla_{i+1}H(x, Sx, S^2x,...,S^{m-1}x).
\end{equation}
Moreover, $H$ could be replaced by a Hamiltonian that is  concave in the first variable and convex in the other variables, though only $m$-cyclically sub-antisymmetric.

\item If the vector fields $u_1, u_2,..., u_{m-1}$ are $m$-cyclically monotone, then (\ref{rep}) holds with $S$ being the identity.  
\end{enumerate}
\end{theorem}

We note first that the general multi-marginal version of the Monge-Kantorovich problem {\bf (K)} where 
the probability measures $\mu_i$, $i=0,1,,...,m-1$ are given marginals on domains $\Omega_i \subset \mathbb{R}^d$, and where $c$ is any  bounded lower semi-continuous cost function $c:\Omega_0 \times \Omega_1\times.... \times \Omega_{m-1}\rightarrow \mathbb{R}$, has the following dual problem.

\newtheorem{duality}{Proposition}[section]
\begin{duality}\label{du}
There exists a solution $\bar \theta$ to the Kantorovich problem (\textbf{K}), as well as an $m$-tuple of functions $(u_0,u_1,..u_{m-1})$ such that for all $i=0, ..., m-1,$
\begin{equation}
u_i(x_i)=\inf_{\substack{x_j \in \Omega_j\\ j\neq i}}\Big( c(x_0,x_1,...,x_{m-1})-\sum_{j \neq i}u_j(x_j)\Big),
\end{equation}
and which maximizes the following dual problem
\begin{equation*}
\sum^{m-1}_{i=0}\int_{\Omega_i} u_i(x_i)d\mu_i \tag{\textbf{D}}
\end{equation*}
among all $m$-tuples $(u_0,u_1,...,u_{m-1})$ of functions $u_i \in L^1(\mu_i)$ for which $\sum_{i=0}^{m-1} u_i(x_i ) \leq c(x_0,...,x_{m-1})$ for all $(x_0,...,x_{m-1}) \in \Omega_0 \times \Omega_1 \times ...\times \Omega_{m-1}$.  

 Furthermore, the maximum value in \textbf{(D)} coincides with the minimum value in \textbf{(K)}, and    \begin{equation}
\hbox{$  \sum\limits_{i=0}^{m-1}u_i(x_i) =c(x_0,...,x_{m-1})$ \quad  for all  $(x_0,...,x_{m-1}) \in support(\bar\theta)$.}
  \end{equation}
\end{duality}
In their seminal paper, Gangbo and \'Swi\c ech \cite{GS} dealt with the case of a quadratic cost function, 
 \begin{equation}\label{quadra}
\hbox{$c(x_0,x_1,x_2,...,x_{m-1}) =\sum_{i=0}^{m-1}\sum_{j=i+1}^{m-1}|x_i-x_j|^2$ on \quad $\mathbb{R}^d \times \mathbb{R}^d \times....\times\mathbb{R}^d$,}
  \end{equation}
and established the following remarkable result.
\begin{theorem} \label{GS} {\rm \bf  (Gangbo-\'Swi\c ech)} Consider Borel probability measures $\mu_i$ on domains $\Omega_i \subset \mathbb{R}^d$, for $i=0,1,...,m-1$ vanishing on $(d-1)$-rectifiable sets and having finite second moments, and let $c$ be a quadratic cost as in (\ref{quadra}). Then,
\begin{enumerate}
\item There exists a unique measure $\bar\theta$ on $\Omega_0 \times \Omega_1 \times ...\times \Omega_{m-1}$, where ({\bf K}) is achieved. It is of the form $\bar\theta=(T_0,T_1,T_2,...,T_{m-1})_{\#}\mu_0$ on $\Omega_0 \times \Omega_1 \times...\times \Omega_{m-1}$, where $T_0=I$, $T_i:\Omega_0 \rightarrow \Omega_i$ and $T_i{}_{\#}\mu_0=\mu_i$.  
\item  
Each  $T_i$ is one-to-one $\mu_i$-almost everywhere, is uniquely determined, and
has the form 
\begin{equation}
\hbox{$T_i(x) = \nabla f_i^*(\nabla f_1(x))$ where
$f_i(x) =|x|^2/2 + \phi_i(x)$ \,  for $x \in \mathbb{R}^d$,}
\end{equation}
and $\phi_i$ is a convex function which is related to the solutions $(u_i)_{i=0}^{m-1}$ of the dual problem  ({\bf D}) by the formula $\phi_i(x)=\frac{m-1}{2}|x|^2 -u_i(x)$ for $x\in \mathbb{R}^d$.
\item Moreover,  $\nabla f_0 (x)=x+T_1x +T_2x +...+T_{m-1}x$ for $\mu_0$-almost all $x \in \mathbb{R}^d$. 
\end{enumerate}
\end{theorem}
The above result was clarified further by Agueh and Carlier \cite{AC}, who essentially established the following. 

\begin{proposition} {\rm \bf  (Agueh-Carlier)} Under the conditions of the Theorem \ref{GS} and with the same notation, we have for each $i=0,..., m-1$, that $\frac{1}{m}\nabla f_i$ is the Brenier map that pushes the measure $\mu_i$ onto the measure $\nu$ on $\mathbb{R}^d$ which is the image of the optimal measure $\bar \theta$ by the ``barycentric map" $(x_0,..., x_{m-1}) \to \frac{1}{m}\sum\limits_{i=0}^{m-1} x_i$. Moreover, $\nu$ is the unique minimizer of the functional $\nu \to \sum_{i=0}^{m-1}W_2^2(\mu_i, \nu)$ where $W_2$ is the Wasserstein distance. 
\end{proposition}

We now describe the situation in the case where the measures $\mu_i$ are obtained from one measure $\mu$ by cyclic permutations.

\begin{corollary} \label{GSS} Let $\mu$ be a probability measure on $\mathbb{R}^N$ and let $\sigma$ be a unitary linear $m$-involution on $\mathbb{R}^N$, that is $\sigma^*=\sigma^{-1}$ and $\sigma^m(x)=x$. Consider the corresponding Kantorovich problem associated to the measures $\mu_i:=\sigma^i_\#\mu$, $i=0,..., m-1$ with quadratic cost on $\mathbb{R}^{mN}$. 
Then,
\begin{enumerate}
\item  The barycentric measure $\nu$ on $\mathbb{R}^N$ associated to the measure $\mu_i:=\sigma^i_\#\mu$, $i=0,..., m-1$ is $\sigma$-invariant. 
\item There exists a strictly convex function  $f_0:\mathbb{R}^N\to \mathbb{R}$ such that the functions $f_i(x):=f_0(\sigma^{m-i}(x))$ satisfy 
\begin{equation}
\hbox{$\nabla f_i{}_\# \mu_i =\nu$ \quad  for  $i=0,..., m-1$.}
\end{equation}
\end{enumerate}
\end{corollary}
\indent{\bf Proof:} As shown by Agueh-Carlier \cite{AC}, $\nu$ is the unique minimizer of the functional $\nu \to \sum_{i=0}^{m-1}W_2^2(\mu_i, \nu)$, where $\mu_i:=\sigma^i_\#\mu$. Since $\sigma^m=I$, the uniqueness yields that $\nu$ is then $\sigma$-invariant. Since now both functions $f_i$ and $\psi_i:=f_0\circ \sigma^{m-i}$ are convex and since both $\nabla f_i$ and $\nabla \psi_i=\sigma^i \circ \nabla f_0\circ \sigma^{m-i}$ push $\mu_i$ onto $\nu$, we get from the uniqueness property of Brenier maps that --modulo a constant-- $f_i=f_0\circ \sigma^{m-i}$.\hfill $\Box$

\subsection*{$m$-cyclically monotone operators}

We shall now apply the above corollary to the case where $\sigma$ is the cyclic permutation $\sigma (x_0, ...., x_{m-1})=(x_1,..., x_{m-1}, x_0)$ on a product space $X^m$. But first, we point to the connection with 
$m$-cyclically monotone operators studied recently by Galichon-Ghoussoub \cite{GG}. They proved the following which is obviously an extension of Krauss' theorem to the case when $m\geq 3$. 

\begin{theorem} (Galichon-Ghoussoub) Let $u_1,..., u_{m-1}:\Omega \to \mathbb{R}^d$ be $m$-cyclically monotone vector fields.  Then, there exists  a Hamiltonian $H\in {\mathcal H}_m$ such that 
\begin{equation}
\hbox{$(u_1(x),..., u_{m-1}(x))=\nabla_{x_2,...,x_m}  H(x,x,...,x)$ for all $x\in \Omega$.}
\end{equation}
Moreover, $H$ can be replaced by a Hamiltonian $K$  on $\mathbb{R}^d\times (\mathbb{R}^d)^{m-1}$, which is concave in the first variable, convex in the last $(m-1)$ variables, whose restriction to  $\Omega^m$ is $m$-cyclically sub-antisymmetric.
 The concave-convex function $K$ is $m$-cyclically antisymmetric in the following sense: For every ${\bf x}=(x_0,..., x_{m-1})$ in $\Omega^m$, we have
\begin{equation}
K(x_0, x_1, ...,x_{m-1})+K_{2,..., m}(x_0, x_1,..., x_{m-1})=0
\end{equation}
where $K_{2,..., m}$ is the concavification of the function $L({\bf x})=\sum\limits_{i=1}^{m-1}H(\sigma^i{\bf x})$ with respect to the last $m-1$ variables. 
\end{theorem}
As to the connection to measure preserving $m$-involutions, they also showed that $u:\Omega \to \mathbb{R}^d$ is $m$-cyclically monotone $\mu$ a.e. if and only if it  is  in the polar set of ${\mathcal S}_m(\Omega, \mu)$, that is 
$
\inf\{\int_\Omega \langle u(x), x-Sx\rangle d\mu; S \in {\mathcal S}_m(\Omega, \mu)\}=0.
$
Equivalently, the projection of $u$ on ${\mathcal S}_m(\Omega, \mu)$ is the identity map, i.e.,  
\[\inf\{\int_\Omega |u(x)-Sx|^2 d\mu (x); S\in {\mathcal S}_m(\Omega, \mu)\}=\int_\Omega |u(x)-x|^2 d\mu (x).
\]
It is easy to see that the above is also equivalent to the statement that 
\begin{equation}
\sup\{\int_{\Omega^m}\langle u(x_0), x_{m-1}\rangle  d\pi({\bf x});\, \pi\in \P_{\rm sym}(\Omega^m, \mu)\}=\int_{\Omega}\langle u(x), x\rangle\, d\mu (x),
\end{equation}
and that the sup is attained at the image of $\mu$ by the map $x\to (x,x,...,x)$, which is nothing but a particular case of the symmetric Monge-Kantorovich problem, when the cost function is given by $c(x_0,..., x_{m-1}) = \langle u(x_0), x_{m-1}\rangle$, where $u$ is an $m$-cyclically monotone operator. 

\subsection*{Multidimensional Monge theorem on graphs and $m$-involutions}

We shall now make a connection between Theorem \ref{GM} above, which is a Monge-Kantorovich problem on $\P_{sym}(\Omega^m,\mu)$ with 
\[
c(x_0, x_1,..., x_{m-1})=\langle u_1(x_0), x_1\rangle +....+\langle u_{m-1}(x_0), x_{m-1}\rangle
\]
as a cost function, and the mass transport result of Gangbo-\'Swi\c{e}ch, which corresponds to the standard multidimensional Monge-Kantorovich problem, i.e., with cost function
\[
c(y_0, y_1,..., y_{m-1})=\sum_{i=0}^{m-1}\sum_{j=i+1}^{m-1}|y_i-y_j|^2, 
\]
and where the marginals are $\mu_i:=\sigma ^i_\#\mu_0$ and $\mu_0$ is supported on an appropriate graph dictated by the vector fields $u_1,..., u_{m-1}$. 

For simplicity, we shall do this for $m=3$, that is for two vector fields $u_1, u_2$  in $L^\infty (\Omega; \mathbb{R}^d)$. Consider the (degenerate) probability measure $\mu$ to be the image of Lebesgue measure $dx$ on $\Omega \subset \mathbb{R}^d$ by the map $P:\Omega \subset \mathbb{R}^d \to \mathbb{R}^{2d} \times \mathbb{R}^{2d} \times  \mathbb{R}^{2d}$ defined by 
$$
x\to    P(x)= (x, x, u_1(x), 0, 0, u_2(x)).
$$
We now consider the $3$-cyclic permutation on $\mathbb{R}^{2d} \times \mathbb{R}^{2d} \times  \mathbb{R}^{2d}$ defined by 
$$
   \sigma \Large((x_{0,1}, x_{0, 2}), ( x_{1,1}, x_{1, 2}), (x_{2,1}, x_{2, 2}))
 = ((x_{1,1}, x_{1, 2}), (x_{2,1}, x_{2, 2}), (x_{0,1}, x_{0, 2})\Large),
$$
in such a way that $\sigma^3 = \Id$.

The quadratic $3$-dimensional Monge problem of Gangbo-\'Swi\c{e}ch applied to the measures $\mu_0=\mu$, $\mu_1=\sigma_\# \mu$ and $\mu_2=\sigma^2_\# \mu$ becomes the problem of minimizing 
$$
   \int_{\Omega}\left\{\|P(x) - \sigma P(S_1x)\|^2
 + \|\sigma P(S_1x) - \sigma^2 P(S_2x)\|^2
 + \|\sigma^2 P(S_2x) - P (x)\|^2\right\}\, dx
$$
over all measurable maps $(S_1,S_2)$, where $S_i: \Omega \rightarrow \Omega$ is measure preserving for $i=1,2$.

The Kantorovich formulation of the problem is then to minimize 
$$
C(\pi):=   \int_{(\mathbb{R}^d)^3}\left\{\|P(x) - \sigma P(y)\|^2
 + \|\sigma P(y) - \sigma^2 P(z)\|^2
 + \|\sigma^2 P(z) - P (x)\|^2\right\}\d \pi(x, y, z)
$$
over all probability measures $\d \pi(x, y, z)$ whose 3 marginals are Lebesgue measure.
In other words,
\begin{eqnarray*}
C(\pi)&=& \int_{(\mathbb{R}^d)^3}\huge\{ |x - u_1(y)|^2 + |x|^2 + |u_1(x)|^2 + |u_2(y)|^2 + |y|^2 + |u_2(x) - y|^2\\
&& \qquad \qquad + |u_1(y)|^2 + |u_2(z)|^2 + |z|^2 + |u_2(y) - z|^2 + |y - u_1(z)|^2 + |y|^2\\
&& \qquad \qquad + |x|^2 + |x - u_2(z)|^2 + |u_1(x) - z|^2 + |z|^2 + |u_1(z)|^2 + |u_2(x)|^2\huge\}\d \pi(x, y, z).
\end{eqnarray*}
Since the integrals of $x^2, y^2, z^2, u_1(x)^2, u_1(y)^2, u_1(z)^2, u_2(x)^2, u_2(y)^2, u_2(z)^2$ against the given marginals of $\pi$ are given constants, the above problem amounts to 
 minimize
 $$
  \int_{(\mathbb{R}^d)^3} \{ |x - u_1(y)|^2 + |u_2(x) - y|^2 
 + |u_2(y) - z|^2 + |y - u_1(z)|^2
 + |x - u_2(z)|^2 + |u_1(x) - z|^2\} \, d\pi (x,y,z),
$$
or, for the same reasons, to maximize
$$
D(\pi):= \int_{(\mathbb{R}^d)^3} \left\{ \langle u_1(y), x\rangle  + \langle u_2(x),  y\rangle  + \langle u_2(y),  z \rangle + \langle u_1(z), y \rangle + \langle u_2(z), x\rangle  + \langle u_1(x), z \rangle \right\} \, d\pi (x,y,z),
$$
which is exactly the problem ${\bf (K_{\rm sym})}$ where the cost 
\[
c(x,y,z)= \langle u_2(x),  y\rangle + \langle u_1(x),  z\rangle
\]
 has been symmetrized. Consider now the optimal maps $ T_1, T_2$ obtained by Gangbo-\'Swi\c{e}ch, that is $ T_i=\nabla f_i^* \circ \nabla f_0$ pushes $\mu_0$ to $\mu_i=\sigma^i_\# \mu_0$ and where $f_i$ is strictly convex for $i=1,2$. By Corollary \ref{GSS}, we have that $f_0=f_i\circ \sigma^i$ for $i=1,2$. 
 
Note now that there exist measure preserving transformations $S_1, S_2$ on $\Omega$ such that  the optimal maps
$$
 \tilde T_1x:= T_1 (Px)=\nabla f_1^*\circ \nabla f_0(Px)= \sigma P (S_1 x)  = (u_1(S_1 x), 0, 0, u_2(S_1 x), S_1 x, S_1 x)
$$
maps $ dx$ onto $\mu_1$, while  
 $$
 \tilde T_2x:= T_2 (Px)=\nabla f_2^*\circ \nabla f_0(Px)
 = \sigma^2\circ P (S_2 x)
 = (0, u_2(S_2x), S_2 x, S_2 x, u_1(S_2 x), 0)
$$
maps $dx$ onto $\mu_2$.

Let now $T_{2,1}:=\nabla f_2^*\circ \nabla f_1$ in such a way that $\tilde T_2=T_{2,1}\circ \tilde T_1$, and let $S_{2,1}$ be a measure preserving on $\Omega$ such that 
\[
T_{2,1}(u_1(y),0,0, u_2(y), y, y)=(0, u_2(S_{2,1}y), S_{2,1}y, S_{2,1}y, u_1(S_{2,1}y), 0), 
\]
so that 
$$T_{2,1}\circ \tilde T_1x= (0, u_2(S_{2,1}\circ S_1 x), S_{2,1}\circ S_1x, S_{2,1}\circ S_1x,  u_1(S_{2,1}\circ S_1 x), 0).
$$
Since  $\tilde T_2=T_{2,1}\circ \tilde T_1$, we have that $S_{2,1}\circ S_1=S_2$, and since $\nabla f_0=\sigma^{3-i} \circ \nabla f_i\circ \sigma^i$ for $i=1, 2$ and $T_{2,1}\circ \sigma (Px)=\nabla f_2^*\circ \nabla f_1\circ \sigma (Px)$ for a.e. $x\in \Omega$,
one can easily verify that 
\[
\hbox{$S_{2,1}=S_1:=S$, $S_2=S_{2,1}\circ S_1=S^2$ and $S^3=I$. }
\]
In other words, the points $(t_0, T_1t_0,  T_2t_0 )$  are such that 
$$
   t_0 
 = P(x) 
 = (x, x, u_1(x), 0, 0, u_2(x)),
$$
$$
  T_1t_0
 = \sigma P (S x)
 = (u_1(S x), 0, 0, u_2(S x), S x, S x)
$$
and
$$
  T_2t_0
 = \sigma^2 P (S^2 x)
 = (0, u_2(S^2 x), S^2 x, S^2 x, u_1(S^2 x), 0)
$$
where $S$ is a measure preserving transformation such that $S^3 = I$.

The convex function $\Phi_0(t)=f_0(t) -\frac{1}{2}|t^2|$  is such that 
$$
   \nabla \Phi_0(t_0) 
 =  T_1t_0 + T_2t_0
 = (u_1(S x), u_2(S^2 x), S^2 x, S^2 x + u_2(S x), S x + u_1(S^2 x), S x).
$$
Define now the convex Lagrangian
$$
 L(x, y, z) = \Phi_0(x, x, y, 0, 0, z) =: \Phi (t).
$$
We then get
$$
 \nabla L(x, y, z) 
 = \bigl( 
    D_{0,1} \Phi(t) + D_{0,2} \Phi(t),
    D_{1, 1} \Phi(t),
    D_{2, 2} \Phi(t) 
   \bigr)
$$
and
$$
   \nabla L (x, u_1(x), u_2(x))
 = (u_1(S x) + u_2(S^2 x), S^2 x, S x).
$$
Let now $H$ be the Legendre transform of $L$ with respect to the last two variables, that is 
\[
H(x,y,z)=\sup\{\langle y, p\rangle +\langle z, q\rangle -L(x, p, q); p\in \mathbb{R}^d, q \in \mathbb{R}^d\}.
\]
$L$ is clearly concave in the first variable, convex in the last two variables and 
\[
L(x, u_1(x), u_2(x)) +H(x, S^2x, Sx) =\langle u_1(x), S^2x\rangle +\langle u_2(x), Sx\rangle,
\]
and in other words, $(u_1(x), u_2(x))=\nabla_{2,3}H(x, S^2x, Sx)$ for all $x\in \Omega$. 

\bibliographystyle{plain}

\end{document}